# Algebraic Properties of Propositional Calculus

Bernd R. Schuh

Dr. Bernd Schuh, Bernhardstraße 165, D-50968 Köln, Germany; bernd.schuh@netcologne.de

*Abstract*.   In this short note we relate some known properties of propositional calculus to purely algebraic considerations of a Boolean algebra. Classes of formulas of propositional calculus are considered as elements of a Boolean algebra. As such they can be represented by uniquely defined elements of this algebra which we call "logical primes". The algebraic notations appear useful because they make it possible to derive well known properties of propositional calculus by simple calculations or to substitute lengthy logical considerations by schematic algebraic manipulations.

*Introduction*

Introductions to the problem of satisfiability can be found in textbooks and reviews, some of them available in the net (see e.g. [1],[2]). One of the unsolved questions of the field is whether satisfiability can be determined in polynomial time ("P=NP ?"). Other questions center around efficient techniques to determine satisfying assignments (see [3,4] for new approaches), and to identify classes of "hard" problems which inherently seem to consume large computing time. I believe that some insight into the difficulties can be gained by using algebraic tools. I will outline these tools in chapters "definitions" and "consequences". Contact wih propositional calculus is made in chapter "propositional calculus". In particular I discuss the representation of formulas in terms of "logical primes" and introduce a group of transformations which leave the number of satisfying assignments invariant.

*Definitions*

We consider a finite algebra V with two operations + and x, and denote by $\mathbf{1}$ and $\mathbf{0}$ their neutral elements, respectively, i.e.

(1)     $a \times \mathbf{1} = a$,  $a + \mathbf{0} = a$

Additionally, the operations are associative and commutative, and the distributive law

(2)     $a \times (b+c) = a \times b + a \times c$

is assumed to hold in V.

Two more properties are required, namely:

(3)     $a + a = \mathbf{0}$

(4)     $a \times a = a$

It is clear from these definitions that V may be identified with the Boolean algebra of propositional calculus, where "x" corresponds to the logical "AND" and "+" to the logical "XOR" (exclusice OR).

To each element of V we introduce its "negation" by

(5)     $\sim a := a + \mathbf{1}$

From (2), (3) and (4) it is clear that $\sim a \times a = \mathbf{0}$ as is appropriate for a negation.

*Consequences.*

As a first consequence of equ.s (1) - (5) we can state the following theorem:

(TI)   $\dim(V) = |V| = 2^N$   for some natural number N

i.e. the number of elements of V is necessarily a power of 2.

This is not surprising, of course, if one has the close resemblence of V to propositional calculus in mind. But here it is to be deduced solely from the algebraic properties.

All proofs are given in the appendix.

In order to formulate a second consequence it is necessary to introduce the notion of "logical primes". We define $p \in V$ as a (logical) prime, iff for any $a \in V$ $pxa = O$ implies $a = O$ or $a = \sim p$. If not clear by definition, the name "prime" will become clear by the following theorem

(TII)   There are exactly $ld|V| = N$ many primes in V. And:

(TIII)  Each element of V has a unique decomposition into primes:

(6)   $a = \Pi_j p_j$   where the product refers to x, and $j \in I_a$, and $I_a = I_b$ iff $a = b$

This property can be formulated alternatively with the negated primes $\sim p_j$ via

(7)   $a = \Sigma_j \sim p_j$ with $j \in {}^c I_a$ (${}^c I_a$ is the complement of $I_a$ in $\{0, 1, ..., N-1\}$ )

The neutral elements *0* and *1* are special cases. *1* is expressed as the empty product according to (6), whereas the sum extends over all primes. For *0* the sum-representation is empty, but the product extends over all possible primes.

A property which is extremely helpful in calculations is

(8)   $\sim p_j \times \sim p_k = \sim p_k \, \delta_{jk}$   ($\delta_{jk} = 1$ iff $j=k$, 0 otherwise)

which with the aid of (5) can be written

$p_j \times p_k = p_j + \sim p_k = \sim p_j + p_k$   for $k \neq j$

Note, that no use has been made of the correspondence of $\{V,+,\times, 0, 1\}$ to propositional calculus, up to now. We can even proceed further and define the analogue of truth assignments. Consider the set of maps $T: V \rightarrow \{0,1\}$. We call T "allowed" iff there is a relationship between the image of a "sum" or a "product" and the image of the single summands or factors. In formula:

(9)   $T(a+b) = f(T(a),T(b))$  and  $T(a \times b) = g(T(a),T(b))$

with some functions f and g and all $a,b \in V$.

These relations suffice to show theorem IV

(TIV)  There are exactly N different allowed maps $T_j$, and they fulfill:

(10)   $T_j(\sim p_k) = \delta_{jk}$

Given functions f and g of (9) one can also use (10) as a definition and extend $T_j$ to all elements of V via (7).

In one last step we assume $N=2^n$ for some natural number n. Then

(TV)  n distinct elements $a_k$ ( different from $0, 1$) can be found, such that

(11)  $\sim p_s = (\Pi_j s_j a_j)(\Pi_k (1-s_k) \sim a_k)$    where $s = \Sigma_r 2^{r-1} s_r$ is the binary representation of s.

In words: each element of V can be written as a "sum" of "products" of all $a_k$ and $\sim a_k$. E.g. for n=3 one has $p_2 = a_2 x \sim a_1 x \sim a_3$ as one of the eight primes. The $a_k$ are not necessarily unique. E.g., for n=3, given $a_k$, the set $a_1, a_3, a_1 x \sim a_2 + \sim a_1 x a_2$ will serve the same purpose (with a different numbering convention in (11)).

*Propositional calculus.*

Propositional calculus (PC) consists of infinitely many formulas which can be constructed from basic variables $a_k$ with logical functions (like "AND", "OR" and negation). Even for a finite set of n basic variables $B_n = \{a_1, a_2, ... a_n\}$ there are infinitely many formulas arizing from combinations of the basic variables. These formulas can be grouped into classes of logically equivalent formulas. That is, formulas F and F' belong to the same class iff their values under any truth assignment $\tau : B_n \rightarrow \{0,1\}$ are the same. Members of different classes are logically inequivalent, i.e. there is at least one truth assignment for which their values differ. This finite set of classes for fixed n can be identified with the algebra V of the foregoing section. Neutral elements of the operations x and +, $1$ and $0$, are interpreted as complete truth and complete unsatisfiability.

In order to see how operations + and x correspond to logical operations "AND" and "OR" we define a new operation v in V via

(12)   a v b = a + b + axb

With this definition the defining relations (1) - (5) can be reformulated in terms of v and x, and the algebraic structure of a Boolean algebra for formulas becomes obvious. v is the logical "OR", x the logical "AND".

Relation (12) reduces logical considerations to simple algebraic manipulations in which + and x can be used as in multiplication and addition of numbers, and additionally the simplifying relations $a+a=\mathbf{0}$ = $\sim axa$ and $axa=a$, $a+\sim a=\mathbf{1}$ hold.

Consider for illustration the so called "resolution" method. It states that avb and $\sim$avc imply bvc. A "calculational" proof of this statement might run as follows (we skip the x-symbol for multiplication in the following and use that in PC the implication $a \Rightarrow b$ is identical to $\sim a$ v b ):

$(avb)(\sim avc) \Rightarrow bvc = \sim((avb)(\sim avc))vbvc =$

$(\mathbf{1}+(a+\sim ab)(\sim a+ca))+b+c+bc+(b+c+bc)(\mathbf{1}+(a+\sim ab)(\sim a+ca)) = \mathbf{1} +ac+\sim ab +$

$+(b+c+bc)(a+\sim ab)(\sim a+ca) = \mathbf{1}+ac+\sim ab+abc+\sim ab+\sim bac = \mathbf{1}+ac(\mathbf{1}+b+\sim b) = \mathbf{1}$

In other words: the implication is a tautology ( true under all truth assignments) as claimed.

TIII and TV tell us that each formula F of PC has a unique decomposition into a "sum" of "products" of its independent variables $a_k$. Because of (8) and (12) the sum in (7) may be written as a "v"-sum. Thus (8) takes the form of a disjunctive normal form (DNF) and it can as well be transformed into a conjunctive normal form (CNF) as given by (6). For the neutral element $\mathbf{0}$ one has

(13)     $\mathbf{0} = (a_1va_2v...va_n)x(\sim a_1v...a_n)x...x(\sim a_1v...\sim a_n)$

with all possible primes. According to (6) each formula F has a similar representation, but with some prime factors missing. From the primes present one can immediately read off the truth assignments for which F evaluates to 0, thus the missing factors give the truth assignments for which F is satisfiable.

Note, however, that each factor in the prime representation of a formula involves *all* $a_k$. So one way of determining satisfying assignments or test a formula for satisfiability consists of transforming a given CNF representation of the formula to its standard form (6). This can be done e.g. by "blowing up" each factor until all $a_k$ are present. E.g.   avbv~c = (avbv~cvd)(avbv~cv~d)  from 3 to 4 variables. Since each new factor has to be treated in the same way, until n is reached, this is a $O(2^n)$ - process in principle, which makes the difficulty in finding a polynomial time algorithm for testing satisfiability understandable.

Also from (7) with (10) and (8) it follows that the satisfying assignments of a formula F= $\Sigma_j$ ~$p_j$ are given by the negated primes which do not show up in the CNF representation. In particular, the number of satisfying assignments is equal to the number of summands in this equation. Furthermore, they can be read off immediately, since, according to (10)  $T_s(F) = 1$ iff the corresponding ~$p_s$ shows up in the sum. Also the Ts must coincide with the $2^n$ possible truth assignments $\tau$:$B_n \rightarrow$ {0,1}. One may choose the numbering such that the values of $T_s$ on $B_n$ are given by the binary representation   s= $\Sigma_r 2^{r-1} T_s(a_r)$.

As a last example for the usefulness of the algebraic approach we consider the number of satisfying assignments of a formula F of PC , #(F) and show that this number does not change if some (or all) of the variables $a_k$ are "flipped", i.e. substituted by their negation and vice versa:

(14)    $\#(F(a_1,...,a_n)) = \#(F(a_1,...\sim a_i,...\sim a_j,...))$

To prove this "conservation of satisfiability" we consider a group of transformations $\{R_0,...R_{N-1}\}$ which negate the $a_k$ according to the following definition: $R_s$ negates all $a_r$ (and $\sim a_r$ likewise) for which $s_r$ in the binary representation of s is non zero. In formula, for any truth assignment $T_j$

$T_j(R_s(a_r)) = (1-s_r) T_j(a_r) + s_r(1- T_j(a_r))$    and    $s= \Sigma_r 2^{r-1} s_r$.

It is easy to see that the $R_s$ form a group with $R_0 = $ id , and each $R_s$ induces a permutation $\pi_s$ of of the $\sim p_j$ which is actually a transposition given by

$\pi_s(j) = s + j - 2\Sigma_r 2^{r-1} s_r j_r$

Thus $R_s$ simply permutes the primes $p_k$ and therefore in the representation of F in (6) or (7) their number is not changed. The fact may also be stated as

$T_j(R_s(F)) = T_{\pi_s(j)}(F)$,

and therefore $\#(F)= \Sigma_j T_j(F) = \Sigma_j T_j(R_s(F)) = \#(R_s(F))$ which proves (14).

*Appendix*

The proofs for theorems (TI) to (TV) are straightforward and only basic ideas will be sketched here.

Proof of TI: For N=1 V consists only of the trivial elements $0$ and $1$. Thus we assume $|V|>2$. For some nontrivial s define $K_s=\{a|axs=0\}$. Obviously $\sim s$ and $0 \in K_s$.

Analogously for $K_{\sim s}$. It is easy to show that $K_s$ and $K_{\sim s}$ are subgroups of V with respect to + , and both have only $O$ in common. Thus each a ε V has a unique decomposition

a=u+v where u ε $K_s$ and v ε $K_{\sim s}$ . Let | $K_s$ |=$N_s$, and | $K_{\sim s}$ |= $N_{\sim s}$. Next we count elements which do not belong to $K_s$ or $K_{\sim s}$. Define:

$E_{Ks}(u_0)$ = {$u_0$+v| v ε $K_{\sim s}$\ $O$} with $u_0$ ε $K_s$. | $E_{Ks}(u_0)$ | = $N_{\sim s}$-1 from the definition. Next

one shows that $E_{Ks}(a)$ and $E_{Ks}(b)$ have no elements in common unless a=b. Thus

|V|= $N_s$-1+ $N_{\sim s}$ +|$\Sigma_u$ $E_{Ks}(u)$ |= $N_s$-1+ $N_{\sim s}$ +( $N_s$-1)|$E_{Ks}(u)$|= ( $N_s$-1)(1+ $N_{\sim s}$-1)+ $N_{\sim s}$

= $N_s$ $N_{\sim s}$ .

Since both $K_s$ and $K_{\sim s}$ are subfields of V (with neutral elements $\sim s$ and s with respect to x) one can apply the same line of argument to each of them until one reaches the trivial field $V_0$={ $O$, $1$} which has |$V_0$|=2. Thus both $N_s$ and $N_{\sim s}$ , and therefore |V| is a power of 2.

Next the proof of (TII) can proceed via induction over N=ld(|V|).

Again one considers the subfields $K_s$ and $K_{\sim s}$ of a V with | V|= $2^{N+1}$ and their sets of primes $p_j$ and $q_j$ which exist by assumption. Then one shows that all $p_j$ + s are primes in V, and $q_j$ + $\sim$s dto. Furthermore one can show that no two of these primes of V or their negations coincide, and, secondly, that any possible prime of V is necessarily one of them. Thus the $p_j$ + s and $q_j$ + $\sim$s constitute the set of primes of V, and their number is by assumption ld($N_s$)+ld($N_{\sim s}$) = N+1.

The fact that different negated $p_k$ are orthogonal, equ. (8), is proven as follows:

For i≠j $p_j$x$\sim p_i$ ε $K_{pi}$ by definition of K. But since $p_i$ is prime, $K_{pi}$ = { $O$,$\sim p_i$}. Thus either $p_j$x$\sim p_i$ = $O$ which implies (because also $p_j$ is prime) that $\sim p_i$ is either $O$ or equal to ~$p_j$ both in contradiction to assumptions, therefore : or $p_j \times$ ~$p_i$ = ~$p_i$ . Which is equivalent to the claim.

Along the same line of thought - considering $K_s$ and $K_{\sim s}$ for s=some prime element of V - it can be proven that each element of V has a unique decomposition into primes, equ. (7) or (6).

Proof of (TIV).

First note that both functions f(x,y) and g(x,y) in equ. (9) can take values 0 or 1 only, and they are symmetric because of the commutativity of the operations x and +. Then from (1) and (9) setting T(a)=0 or 1 respectively one gets

$$0 = g(0, T(\mathbf{1})) = g(T(\mathbf{1}), 0) \quad \text{and} \quad 1 = g(1, T(\mathbf{1})) = g(T(\mathbf{1}), 1) \quad \text{and}$$

$$T(\mathbf{0}) = g(1, T(\mathbf{0})) = g(0, T(\mathbf{0})) \quad \text{from ax } \mathbf{0} = \mathbf{0}.$$

If one chooses $T(\mathbf{0}) = 0$ then $T(\mathbf{1}) = 0$ leads to a contradiction, as well as setting both values equal to 1. One is left with the choice

(A)    $T(\mathbf{0}) = 0$ and $T(\mathbf{1}) = 1$

(B)    $T(\mathbf{0}) = 1$ and $T(\mathbf{1}) = 0$

We adopt choice (A) in the following. As a consequence

0=g(0,1) = g(1,0) = g(0,0) and 1 = g(1,1) and, from (1) for +

0=f(0,0)=f(1,1) and 1=f(1,0)=f(0,1) .

Let T be fixed. Because of (8): 0=g(T(~p),T(~q)) for different p,q. Thus either T(~p)=T(~q)=0 or the two assignments have different value. If T(~$p_k$)=0 for all k, one gets a contradiction to $\mathbf{1} = \Sigma_k$ ~$p_k$ and 0=f(0,0). Thus at least for one k T(~$p_k$)=1. But then for all other j T(~$p_j$)=0 because of 0=g(0,1) and the orthogonality relation (8).

Thus for each T there is exactly one $\tilde{p}_k$ with truth assignment 1, and all other $\tilde{p}$ giving 0. Now consider two different maps T, T' with $T(\tilde{p}_k)=1$ and $T'(\tilde{p}_l)=1$. Then k and l must be different, otherwise the two maps would coincide. Repeating this argument with a third T" and so on leads to the conclusion that there are exactly as many allowed maps as there are primes. We can label the maps as we would like to, so the most natural choice is equ. (10).

As for theorem V, the easiest way to prove the existence of $n=ld(N)$ $a_k$ is to construct them from the uniquely defined primes:

$a_r = \Sigma_i \Sigma_s \Sigma_l \tilde{p}_i \delta(i,s+2^k l)$

where $\delta$ is the Kronecker $\delta$ and the s and l sums run from $2^{k-1}$ to $2^k-1$ and from 0 to $2^{n-k}-1$ respectively. Constructing them inductively is more instructive because one encounters choices which lead to different sets of $a_k$. The seemingly complicated formula above is obsolete once one uses the binary representation of all quantities which is given by the bijection $F \longleftrightarrow T_{N-1}(F) ... T_i(F) ... T_0(F)$ for any F. In particular the $a_i$ take the simple form:

$a_1 =$     ....1010101010101010

$a_2 =$     ....1100110011001100

$a_3 =$     ....1111000011110000

and so on.